\RequirePackage{fix-cm}
\documentclass[smallextended]{svjour3}       % onecolumn (second format)
\smartqed  % flush right qed marks, e.g. at end of proof
\usepackage{graphicx}

\usepackage{bm}

\usepackage{amssymb, amsmath} 
\usepackage{latexsym} 
\usepackage{graphics} 
\usepackage{url} 
\usepackage{enumerate}
\usepackage[all]{xy}
\usepackage{float}
\usepackage{amsrefs}
\usepackage{color}
\usepackage{tabularx}
\usepackage{ltablex}

\newtheorem{theo}{Theorem}%[section]

\newtheorem{lemm}[theo]{Lemma}

\newtheorem{defi}[theo]{Definition}

\newenvironment{proo}[1][\proofname]{\normalfont{\itshape
#1{:}}\quad\mdseries\ignorespaces}
{{\hspace{10cm}$\Box$}{\vskip\belowdisplayskip}}

\journalname{}
\begin{document}

\title{The relation between symmetries and coincidence and collinearity of polygon centers and centers of multisets of points in the plane %\thanks{Grants or other notes
}
\subtitle{}%Do you have a subtitle?\\ If so, write it here}

\titlerunning{Symmetries and coincidence of centers}        % if too long for running head

%\author{Marta Farr\'e Puiggali \and L. Felipe Prieto-Mart\'inez %etc.
%}

\author{Luis Felipe Prieto-Mart\'inez}


%\institute{Marta Farr\'e Puiggali \at
   %         \textcolor{red}{FALTAFALTAFALTA}
%\email{FALTAFALTA}
 %   \and
\institute{Luis Felipe Prieto-Mart\'inez \at
    Universidad Polit\'ecnica de Madrid (Spain), Departamento de Matem\'atica Aplicada
%Universidad Aut\'onoma de Madrid (Spain), Department of Mathematics \\
              \email{luisfelipe.prieto@upm.es}           %  \\
       }

\date{}

\maketitle

\begin{abstract}  There are several remarkable points, defined for polygons and multisets of points in the plane, called centers (such as the centroid). To make possible their study, there exists  a formal definition for the concept of center in both cases. In this paper, the relation between symmetries  of polygons and multisets of points in the plane and the coincidence and collinearity of their centers is studied. First, a precise statement for the problem is given. Then, it is proved that, given a polygon or a multiset of points in the plane, a given point in the plane is a center  for this object if and only if it belongs to the set of  points fixed by its group of symmetries.

 \keywords{Polygon \and Triangle Center \and Polygon center \and Finite multiset of points in the plane \and Group of Symmetries }
 %\PACS{\textcolor{red}{PACS code1 \and PACS code2 \and more}}
 \subclass{Primary 51M04  \and Secondary 51M15}
\end{abstract}

\section{Introduction}

Associated to the every triangle $P$ there are four famous points known as the \emph{classical centers}. Those are the incenter, the centroid, the circumcenter and the orthocenter. But there are many other remarkable points, associated to $P$, which are also called “centers”. Thousands of them
are known at this moment (see the web site of the \emph{Encyclopedia of Triangle Centers} \cite{K.ETC}).

But, what must satisfy a point to be called “center”? In a series of articles published in the 1990’s  (see, for example, \cite{K.FE, K.CL}) C. Kimberling indicated the importance of giving a formal
 definition of triangle center as a function and not as a ``remarkable'' point obtained with a geometric construction. Recently, the ideas by Kimberling has been generalized for polygons with four or more sides in \cite{FP, PS}.

In this paper we deal, not only with centers of $n$-gons but with $n$-multiset of points in the plane (multiset with $n$ elements). Polygons are multisets of points in the plane endowed with a notion of adjacency between their elements, called vertices.  Polygons (and their centers) are  basic objects in Elementary Geometry and multisets of points in $\mathbb R^N$ (and their centers) are very important in problems in Applied Mathematics (see for instance \cite{P} to see an application of the centroid in the study of tumor growth).

Concerning this theory of centers, in the bibliography we can find (not explicitely stated) the following:

\medskip

\noindent \textbf{Principle:} \emph{As more coincidence and collinearity of centers occur for a given $n$-gon (resp. $n$-multiset of points in the plane), the more regular it is.}

\medskip

\noindent See for example \cite{AS.Q} for the case of quadrilaterals, \cite{FP, PS} for general $n$-gons  and \cite{E,K} for  $n$-multisets of points in $\mathbb R^N$.

Let us consider, as inspiring examples for our study, the following results:

\begin{itemize}

\item For a triangle with not all its vertices collinear, the incenter and the orthocenter  coincide if and only if the triangle is equilateral (see \cite{I}).

\item For a triangle with not all its vertices collinear, the incenter, the centroid and the orthocenter  are collinear if and only if the triangle is isosceles (see \cite{F}).

\item For a set of four different points $\{V_1,V_2,V_3,V_4\}$ in the plane such that any of its elements is in the convex hull of the other three, the centroid   and the Fermat-Torricelli point  coincide if and only $V_1,V_2,V_3,V_4$ correspond to the vertices of a parallelogram (see \cite{AS.Q}, where the results is stated for quadrilaterals).

\item  For a convex quadrilateral, the centroid  and the centroid of the boundary  coincide if and only if it is a paralellogram (see \cite{AS.Q}).

\end{itemize}

The main target of this paper is, in relation to the \emph{Principle} above, to solve the following:

\medskip

\noindent \textbf{Main Problem:} \emph{For a given polygon $P$ (resp. multiset $\widetilde P$), determine the set of points in the plane that can be viewed as centers, according to the formal definition of center as a function (see Definitions  \ref{defi.centermulti} and \ref{defi.centerpoly}).}

\medskip
\noindent  This problem was already explored (but not solved) in \cite{FP, PS} for $n$-gons  and in \cite{E, K} for $n$-multisets of points in $\mathbb R^N$.  We prove that the answer is (exactly) the set of fixed points by the group of symmetries of $P$ (resp. $\widetilde P$) and is given in Theorem \ref{theo.multiset} (for $n$-multisets of points in the plane) and Theorem \ref{theo.polygon} (for $n$-gons).

The fact that the set of centers must be contained in this set of fixed points is easy to verify and the main ideas already appear in the bibliography cited above (anyway, a proof fitting our approach is included here). But to solve  the Main Problem there are two difficulties: (a) {find two $n$-gon centers such that they coincide if and only if the corresponding $n$-gon has rotational symmetry (resp. for $n$-multisets) and (b) {find three $n$-gon centers  such that they are collinear  if and only if the corresponding $n$-gon  has an axis of symmetry (resp. for $n$-multisets).}

The point is that, for general values of $n$, it is not easy (it may be impossible) to find simple examples of centers with these properties.  It will be neccessary to construct artificially two $n$-multiset centers and two $n$-gon centers with this purpose of detecting asymmetry. This centers also apply for non-covex and non-simple $n$-gons, which is one of the strengths of the results herein.

%This centers are not continuous with respect to the more natural topologies for the set of $n$-gons and for the set of $n$-multisets.

These proccesses of detecting asymmetry can be reduced, in turn, to solving a nice  problem (and its weighted version) concerning cyclic convex $n$-gons, which may be of its own interest.

\medskip

\noindent \textbf{First Step Problem:} \emph{Let $c$ be a circle with center $O$. Find a function $\widetilde{\mathfrak B}_n$ that assign to each set of points $\{V_1,\ldots,V_n\}\subset c$ a point in $c\cup\{O\}$  such that (1) $\widetilde{\mathfrak B}_n(\{V_1,\ldots,V_n\})=O$ if and only if there is a rotation fixing the set $\{V_1,\ldots,V_n\}$  and (2) it commutes with similarities, that is, for every  similarity $T$,  $\widetilde{\mathfrak B}_n$ satisfies}
$$\widetilde{\mathfrak B}_n(\{T(V_1),\ldots,T(V_n)\})=T(\widetilde{\mathfrak B}_n(\{V_1,\ldots,V_n\})). $$

The notation and formal definitions (including the one of $n$-multiset and $n$-gon center) required for the rest of the paper are introduced in Section \ref{section.notation}. The First Step Problem is solved in Section \ref{section.auxiliar}. Finally, in Sections \ref{section.multiset} and \ref{section.polygon} we include the proof of Theorems \ref{theo.multiset} and \ref{theo.polygon} respectively.

\section{Notation and definition of center} \label{section.notation}
%%%%%%%%%%%%%%%%%%%%%%%%%%%%%%%%%%%%%%%%%%%%%%%%%%%%%%%%%%%%%%%%%%%%%

The most basic objets through the rest of this paper are $n$-multiset (multisets of $n$ elements for $n\geq 1$) of points in the plane (just referred as $n$-multiset in the following). The elements will be referred as \textbf{vertices}. If all the elements are different we call  $n$-set to the correspoding object. Let $\widetilde{\mathcal P}_n$ be the set consisting in all $n$-multisets of points in $\mathbb R^2$.

 Eventually (in the next section) we will also denote by $\widetilde{\mathcal P}^*$ to the set of $n$-multisets with their vertices labelled with a natural number, that is, the sets of pairs $\{(V_1,l_1),\ldots,(V_n,l_n)\}$ where $(V_i,l_i)\subset \mathbb R^2\times \mathbb N$.

 For $n\geq 3$, a \textbf{$n$-gon} is a $n$-multiset $\widetilde P$ which elements are called \textbf{vertices} endowed with some notion of \textbf{adjacency} $\mathcal A$ between these points (symmetric and irreflexive) that additionally satisfies (1) each vertex is adjacent to (exactly) two vertices and (2) for every pair of vertices $V,V'$ there exist two (and only two) sequences of vertices $V_1,\ldots,V_k$, that only intersect in their endpoints, such that $V=V_1$, $V'=V_k$ and such that, for $i=1,\ldots,k-1$, $V_i$ is adjacent to $V_{i+1}$. The set of all $n$-gons will be denoted by $\mathcal P_n$.

The \textbf{sides} of a $n$-gon are the segments joining two adjacent vertices.  The rest of segments joining vertices of the $n$-gon are called \textbf{diagonals}. Note that, for the case $n=3$,  every 3-multiset has a unique triangle (3-gon) structure.

In the notation above, any sequence of vertices $(V_1,\ldots,V_k)$ in $\widetilde P$ such that, for $i=1,\ldots,k-1$, $V_{i+1}$ is adjacent to $V_i$  is called a \textbf{polygonal chain} starting at $V_1$ (sometimes, this concept refers to the set of sides joining these vertices). This polygonal chain is \textbf{closed} if $k=n+1$ and so $V_1=V_k$. Let $r_1,\ldots, r_{k-1}$ where $r_i$ denotes the side with endpoints $V_i,V_{i+1}$. We say that the polygonal chain is \textbf{simple} if these segments only intersect in their endpoints and only with the following and the previous one. For the case $k=n+1$ we consider $r_1$ to be ``the following one'' of $r_n$ and, if the chain is simple, we say that the polygon $P$ is \textbf{simple}.

% and is \textbf{simple} if so is the union of the corresponding sides (as a curve).

\medskip

\noindent \textbf{Comment on the definition of polygon:} \emph{We may assume that the $n$-gon has a labelling for its vertices $(V_1,\ldots,V_n)$ in such a way that $V_i,V_j$ are adjacent if and only if $i,j$ are consecutive (modulo $n$). This labelling is not unique. Let us denote by $S_n$ to the set of {permutations} of $\{1,\ldots, n\}$. Consider $D_n<S_n$ to be the {dihedral subgroup}, that is, the one generated by the permutations $\rho,\sigma$ given by}
$$\rho(i)=i+1\mod n,\qquad\sigma(i)=n-i+2\mod n.$$

\noindent \emph{Then two labellings correspond to the same $n$-multiset if and only if there exists some $\alpha\in D_n$ or $\alpha\in S_n$ in each case, by}
$$(V_1,\ldots, V_n)\longmapsto (V_{\alpha(1)},\ldots,V_{\alpha(n)}) $$

\noindent \emph{This second approach is more similar  to the one in \cite{ES, FP, PS}.}

\medskip

To avoid confusion, in general, the objects related to $n$-multisets will be denoted with tildes and the ones related to $n$-gons without them.

Denote by $E(2)$, $S(2)$, with $E(2)<S(2)$, to the groups of plane rigid motions and plane similarities, respectively.  $E(2)$, $S(2)$ act on $\widetilde{\mathcal P}_n$, $\mathcal P_n$ and $\widetilde{\mathcal P}_n^*$. The action of any $T\in S(2)$ is given by
$$\begin{array}{l l l}\{V_1,\ldots,V_n\}&\longmapsto& \{T(V_1),\ldots T(V_n)\} \\
(\{(V_1,\ldots,V_n\},\mathcal A)&\longmapsto& (\{T(V_1),\ldots T(V_n)\},T(\mathcal A)) \\
&&\text{ where }(W_1,W_2)\in T(\mathcal A)\Leftrightarrow (W_1,W_2)\in\mathcal A \\
\{(V_1,l_1),\ldots,(V_n,l_n)\}&\longmapsto& \{(T(V_1),l_1),\ldots (T(V_n),l_n)\} \end{array}$$

%\noindent  We say that $\widetilde P,\widetilde Q\in\widetilde{\mathcal P_n}$ (resp. $P,Q\in\mathcal P_n$ or $\widetilde P_n^*,\widetilde Q_n^*\in\widetilde{\mathcal P}_n^*$) are \textbf{congruent} or \textbf{similar} if they both lie in the same orbit with respect to this action with respect to $E(2),S(2)$, respectively.

 The \textbf{group of symmetries} of a given $n$-multiset $\widetilde P$ (resp. of a a $n$-gon $P$ or of a $n$-multiset with its vertices labelled $\widetilde{P}_n^*$) is the subgroup of plane rigid motions $T\in E(2)$ such that $T(\widetilde P)=\widetilde P$ (resp. $T(P)=P$, $T(\widetilde P^*)=\widetilde P^*$).

We say that some $\widetilde P\in\widetilde{\mathcal P}_n$ (resp. $P\in\mathcal P_n$ or $\widetilde{\mathcal P}_n^*$) is \textbf{rotationally symmetric} or order $k$ if there is a rotation of order $k$ in its group of symmetries.   We say that some $\widetilde P\in\widetilde{\mathcal P}_n$ (resp. $P\in\mathcal P_n$ or $\widetilde{\mathcal P}_n^*$) is \textbf{axially symmetric} if there is some reflection with respect to a line in its group of symmetries. 

\medskip

\noindent \textbf{Remark:} \emph{Let $P\in\mathcal P$ and $\widetilde{P}$ be the $n$-multiset associated to $P$ (containing its vertices). Note that the group of symmetries of $P$ is a subgroup of the group of symmetries of $\widetilde P$ and both groups may not be equal. As a consequence, the set of fixed points of the group of symmetries of $P$ contains the set of fixed points by the group of symmetries of  $\widetilde{P}$ and may be larger.}

\medskip

With this notation, we have the following:

\begin{defi} \label{defi.centermulti} For $n\geq 1$, let $\widetilde{\mathcal F}_n$ be a non-empty subset of $\widetilde{\mathcal P}_n$ closed with respect to similarities.  A \textbf{$n$-gon center} is a function $\widetilde{\mathfrak X}_n:\widetilde{\mathcal F}_n\to\mathbb R^2$ that commutes with respect to similarities.

\end{defi}

\begin{defi} \label{defi.centerpoly} For $n\geq 3$, let $\mathcal F_n$ be a non-empty subset of $\mathcal P_n$ closed with respect to similarities.  A \textbf{$n$-gon center} is a function $\mathfrak X_n:\mathcal F_n\to\mathbb R^2$ that commutes with respect to similarities.

\end{defi}

The most important multiset center is the \textbf{centroid} $\widetilde{\mathfrak C}_n:\widetilde{\mathcal P}_n\to\mathbb R^2$, which is given by
$$\widetilde{\mathfrak C}_n(\{V_1,\ldots, V_n\})=\frac{1}{n}V_1+\ldots+\frac{1}{n}V_n. $$

\noindent One more illustrative example (the circumcenter) will be explained at the beginning of the following section.

Associated to every $n$-multiset center $\widetilde{\mathfrak Z}_n:\widetilde{\mathcal F}_n\to\mathbb R^2$ there is  a $n$-gon center $\mathfrak Z_n:\mathcal F_n\to\mathbb R^2$, defined for the family $\mathcal F_n$ of $n$-gons $P$ which $n$-multiset of vertices $\widetilde P$ is in $\widetilde{\mathcal F}_n$ and given by $\mathfrak Z_n(P)=\widetilde{\mathfrak Z}(\widetilde P)$. So we also have a $n$-gon version of the centroid $\mathfrak C_n:\mathcal P_n\to\mathbb R^2$ that maps each $n$-gon to the centroid of its vertices.

The converse is not true: not for very $n$-gon center there is associated some $n$-multiset center. For example, for $n=4$,  the function $\mathfrak D_4:\mathcal F_4\to\mathbb R^2$ that maps, when defined, each tetragon $P$ to its diagonal crosspoint is a 4-gon center but there is no finite multiset center associated to it, since its definition depends on the adjacency relation  $\mathcal A$ of the vertices of the polygon.

\medskip

\noindent \textbf{Comment:} \emph{In this paper, we commit an abuse of notation which is very extended in this context. We call ``center''  to both, the function  $\widetilde{\mathfrak Z}_n$ (resp. $\mathfrak Z_n$) and the point  $\widetilde{\mathfrak Z}_n(\widetilde P)$ corresponding to a given multiset $\widetilde P$ (resp. the point $\mathfrak Z_n(P)$ corresponding to a given  $n$-gon $P$) in its domain.}

%%%%%%%%%%%%%%%%%%%%%%%%%%%%%%%%%%%%%%%%%%%%%%%%%%%%%%%%%%%%%%%
\section{Some centers for families of cyclic $n$-multisets} \label{section.auxiliar}

We say that a $n$-set $\widetilde P$ is \textbf{cyclic} if all its elements belong to a circle. The center of this circle is called the \textbf{circumcenter} of $\widetilde P$. The function $\mathfrak M_n:\mathcal F_n\to\mathbb R^2$, for $\mathcal F_n$ being the family of cyclic $n$-sets, that maps every cyclic $n$-set to its circumcircle is a $n$-multiset center.

Any cyclic $n$-set $\widetilde P$ has a natural notion of adjacency between its elements (so cyclic $n$-multisets endowed with this adjacency relation can be viewed as $n$-gons). For $n\geq 3$, two points $V,W$ in $\widetilde P$ are adjacent if there is no other point of $P$ for some  of the two circular arcs with endpoints $V,W$.

Let us start with the following lemma. It is required to guarantee that some of the algorithms below produce an output.

\begin{lemm} \label{lemm.aux} Let $\widetilde P, \widetilde Q$ be two sets of points inscribed in the same circle $c$ with center $O$ and consisting, respectively, in the vertices of a regular $n$-gon and a regular $m$-gon.

\begin{itemize}

\item[(i)] The set $\widetilde P\cup \widetilde Q$ is not rotationally symmetric if and only if $m,n$ are  coprimes.
\end{itemize}

\noindent Provided that $m,n$ are coprime, then:

\begin{itemize}

\item[(ii)] $\widetilde P,\widetilde Q$ intersect in at most one point.

\item[(iii)] If there are two pairs of points $(V_1,W_1)$, $(V_2,W_2)\in \widetilde P\times \widetilde Q$ such that they are adjacent in $\widetilde P\cup\widetilde Q$ and $\angle (V_1,O,W_1)=\angle(V_2,O,W_2)=\alpha<\pi$, then $\widetilde P\cup \widetilde Q$ is axially symmetric with respect to the segment bisector of $V_1,V_2$ (or, equivalently, of $W_1,W_2$). Moreover, the midpoints $M_1,M_2$ of the circular arcs corresponding to these angles cannot be antipodal. See the right hand side of Figure \ref{lemmaenunciado}.

\item[(iv)] There are not three pairs of adjacent points $(V_1,W_1)$, $(V_2,W_2)$, $(V_3,W_3)\in \widetilde P\times \widetilde Q$ such that $\angle(V_1,O,W_1)=\angle(V_2,O,W_2)=\angle(V_3,O,W_3)$.

\end{itemize}

\end{lemm}

\begin{figure}[h]
\centering
\includegraphics[width=0.85\textwidth]{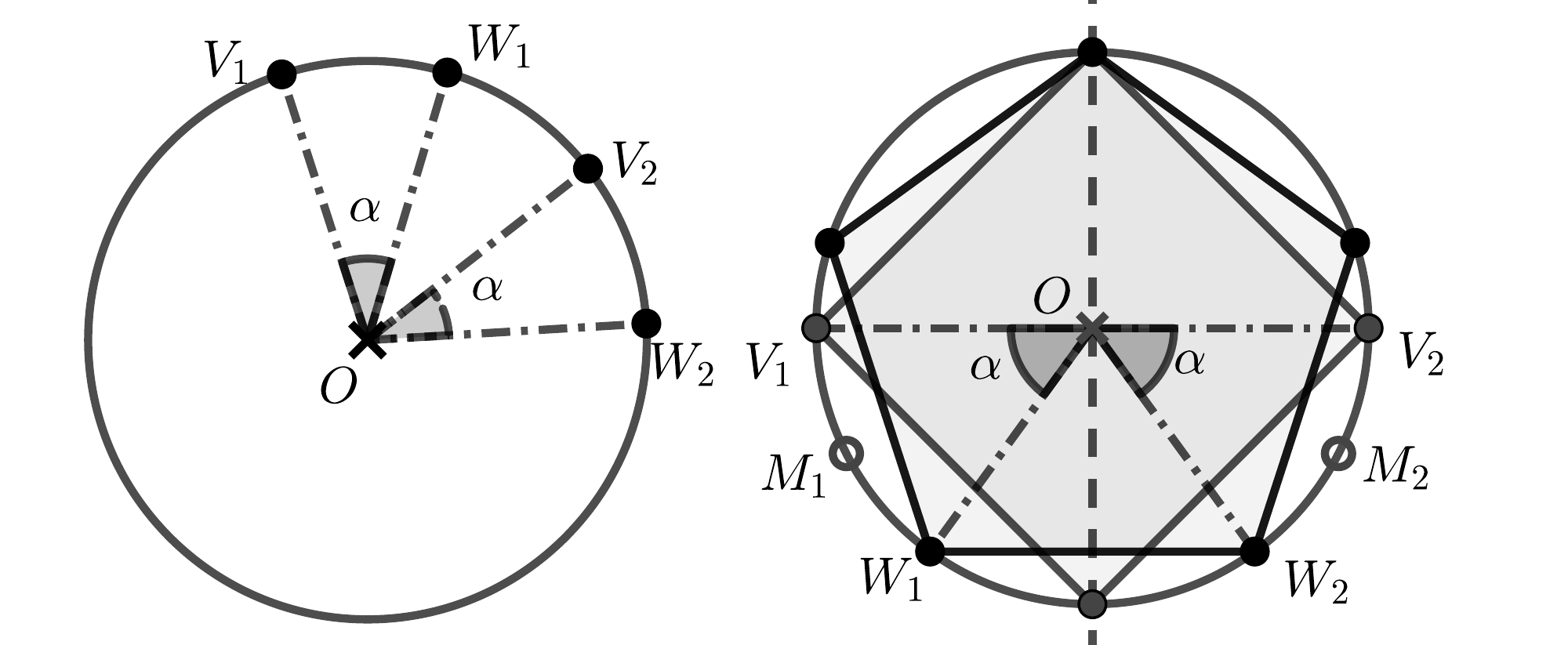}
\caption{Illustrations for Statement (iii). The configuration on the left hand side is not possible.}
\label{lemmaenunciado}
\end{figure}

\begin{proo}  Statements (i) and (ii) are immediate.

\begin{itemize}

\item[(iii)] Suppose that there are two pairs of points $(V_1,W_1),(V_2,W_2)\in \widetilde P\times \widetilde Q$ such that $\angle(V_1,O,W_1)=\angle(V_2,O,W_2)=\alpha$.

Note that $W_1,W_2$ cannot be ``at the same side'' of $V_1,V_2$, respectively (the circular arcs from $V_1$ to $W_1$ and from $V_2$ to $W_2$ cannot have the same orientation, this would correspond to the picture at the left hand side of Figure \ref{lemmaenunciado}). To see this, let $\widetilde P'$  be the set obtained from $\widetilde P$ by a rotation of angle $\alpha$ in the sense of the circular arc from $V_1$ to $W_1$. Then  $\widetilde P',\widetilde Q$ intersect at two points, which is a contradiction with Statement (ii).

So the segment bisector of the segments $V_1,V_2$ and $W_1,W_2$ coincide. To see that $\widetilde P\cup \widetilde Q$ is symmetric with respect to this line, we can use that  any regular $n$-gon is symmetric with respect to any of the segment bisectors of its sides and diagonals.

Finally, suppose that the midpoints of these arcs are antipodal. Then, we are (modulo congruences) in the situation in Figure \ref{lemma}. The angle $\pi-\alpha$ and $\pi+\alpha$ are multiples of $\frac{2\pi}{m}$ and of $\frac{2\pi}{n}$.  So $\alpha$ is a multiple of $\frac{\pi}{m}$ and $\frac{\pi}{n}$. The case $\alpha=0$ is not possible (contradiction with Statement (i)). The case $\alpha=\frac{\pi}{m}=\frac{\pi}{n}$ is a contradiction with the fact that $m,n$ are coprime. Finally, the case $\alpha=k\frac{\pi}{m}$ (resp. $\alpha=k\frac{\pi}{n}$) for $k\geq 2$ contradicts that $V_1,V_2$ are adjacent.

\begin{figure}[h!]
\centering
\includegraphics[width=0.5\textwidth]{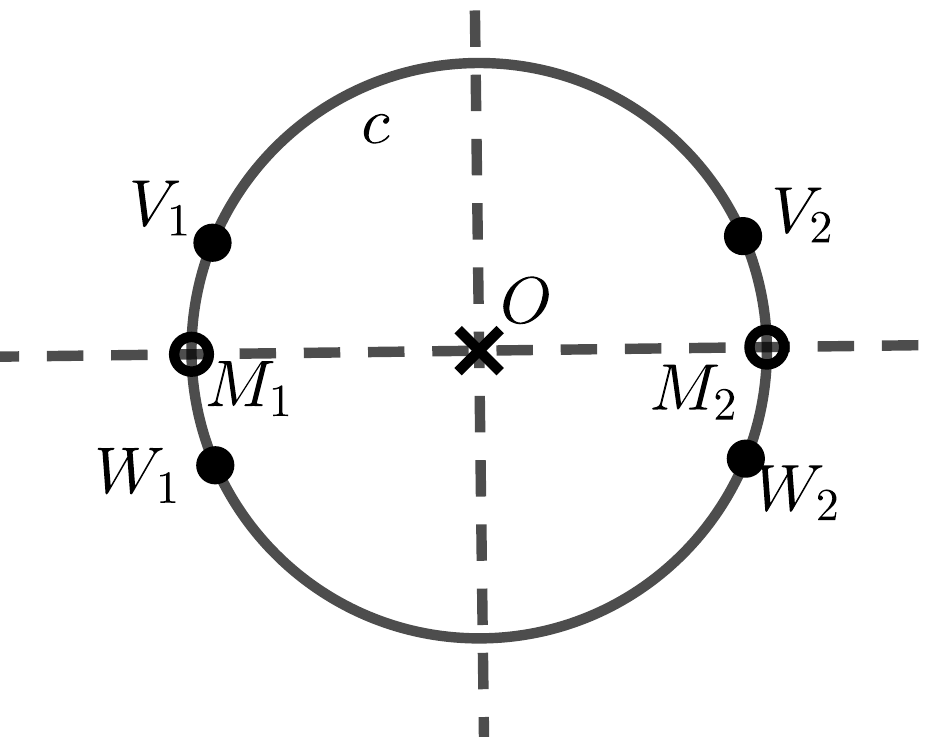}
\caption{$M_1,M_2$ cannot be antipodal.}
\label{lemma}
\end{figure}

%So $2\pi=k_1(\pi-\alpha)m=k_2(\pi-\alpha)n$ and $2\pi=r_1(\pi-\alpha)m=r_2(\pi-\alpha)n$. This is a contradiction with the fact that $m,n$ are coprime.
%\noindent Suppose that there are $k$ points in $\widetilde P$ in the circular arc with endoints $V_1,V_2$ and contained in the upper half-circle. On the other hand, there are $k$ or less  points in $\widetilde P$ in the circular arc with endoints $W_1,W_2$ and contained in the lower half-circle

\item[(iv)] Suppose that there are such three pairs of points. Then two of the points $W_{i_1},W_{i_2}$ for $i_1,i_2\in\{1,2,3\}$ are ``at the same side'' of the corresponding points $V_{i_1},V_{i_2}\in\widetilde P$ in the sense precised before. So we can repeat the argument in the first paragraph of the proof of Statement (iii).

\end{itemize}

\end{proo}

Now we describe our algorithms. Each algorithm describes one of the functions $\Phi$, $\widetilde{\mathfrak A}_n$, $\widetilde{\mathfrak B}_n$, $\widetilde{\mathfrak{B}}_n^*$ which domain is in $\widetilde{\mathcal P}_n$ or $\widetilde{\mathcal P}_n^*$ and which image is in $\widetilde{\mathcal P}_n$ or $\mathbb R^2$. For us the following fact of immediate proof for each case will be very important:

\medskip

\noindent \textbf{Remark:} \emph{The functions $\Phi$, $\widetilde{\mathfrak A}_n$, $\widetilde{\mathfrak B}_n$, $\widetilde{\mathfrak{B}}_n^*$, described in Algorithms 1, 2, 3, 4, commute with respect to similarities (in fact, $\widetilde{\mathfrak A}_n$ and $\widetilde{\mathfrak B}_n$ are $n$-multiset centers). In particular, $\widetilde{\mathfrak B}_n$ is the answer for the \emph{First Step Problem}.}

\noindent \begin{tabularx}{\textwidth}{|l X|} \hline  & \textbf{Algorithm 1 (definition of $\Phi$)} \\

\hline

INPUT:& a cyclic $n$-set   $\widetilde P\in\widetilde{\mathcal P}_n$ contained in a circle $c$ with circumcenter $O$ which is rotationally symmetric of order $k\geq 2$.\\

OUTPUT:& a $k$-set $\Phi(\widetilde P)$ consisting in the vertices of some regular $k$-gon inscribed in $c$ (or the endpoints of a diameter, for $k=2$).\\

\hline (1) & If $\widetilde P$ already is the $k$-set consisting in the vertices of some regular $k$-gon inscribed in $c$, then $\Phi(\widetilde P)=\widetilde P$. \\

 (2) & In other case, for each element $V$ in $\widetilde P$, consider the sequence $V_1,\ldots,V_n$ where $V_1=V$, and, for $i=1,\ldots,n-1$,  $V_{i+1}$ is the element in $\widetilde P$ adjacent to  $V_i$ and such that the circular arc from $V_i$ to $V_{i+1}$ is positively oriented.\\

(3) & For each element $V$ in $\widetilde P$, and for the notation described above, define the sequence $(\alpha_1,\ldots,\alpha_n)$ such that $\alpha_i$ is the angle $\angle(V_i,O,V_{i+1})$ (corresponding to the circular arc from $V_1$ to $V_{i+1}$).\\

(4) & There are exactly $k$ elements in $\widetilde P$ such that the corresponding sequence $(\alpha_1,\ldots,\alpha_n)$ is minimal with respect to the lexicographic order. Let $\widetilde Q_+$ be this $k$-set.\\

(5) & For each element $V$ in $\widetilde P$, consider the sequence $W_1,\ldots,W_n$ where $W_1=V$, and, for $i=1,\ldots,n-1$,  $W_{i+1}$ is the element in $\widetilde P$ adjacent to  $W_i$ and such that the circular arc from $W_i$ to $W_{i+1}$ is negatively oriented.\\

(6) & For each element $V$ in $\widetilde P$, and for the notation described above, define the sequence $(\beta_1,\ldots,\beta_n)$ such that $\beta_i$ is the angle $\angle(W_i,O,W_{i+1})$ (corresponding to the circular arc from $W_i$ to $W_{i+1}$).\\

(7) & There are exactly $k$ elements in $\widetilde P$ such that the corresponding sequence $(\beta_1,\ldots,\beta_n)$ is minimal with respect to the lexicographic order. Let $\widetilde Q_-$ be this $k$-set.\\

(8) & The case $\widetilde Q_+=\widetilde Q_-$ is not possible. %then $\Phi(\widetilde P)=\widetilde Q_+$.\\
So, there are $k$ circular arcs oriented in the positive order starting in an element $V$ in $\widetilde Q_+$ and ending in an element $W$ in $\widetilde Q_-$ adjacent (in $\widetilde Q_+\cup\widetilde Q_-$) to $V$. $\Phi(\widetilde P)$ is the set of midpoints of these circular arcs.\\

\hline

\end{tabularx}

\begin{figure}[h]
\centering
\includegraphics[width=0.7\textwidth]{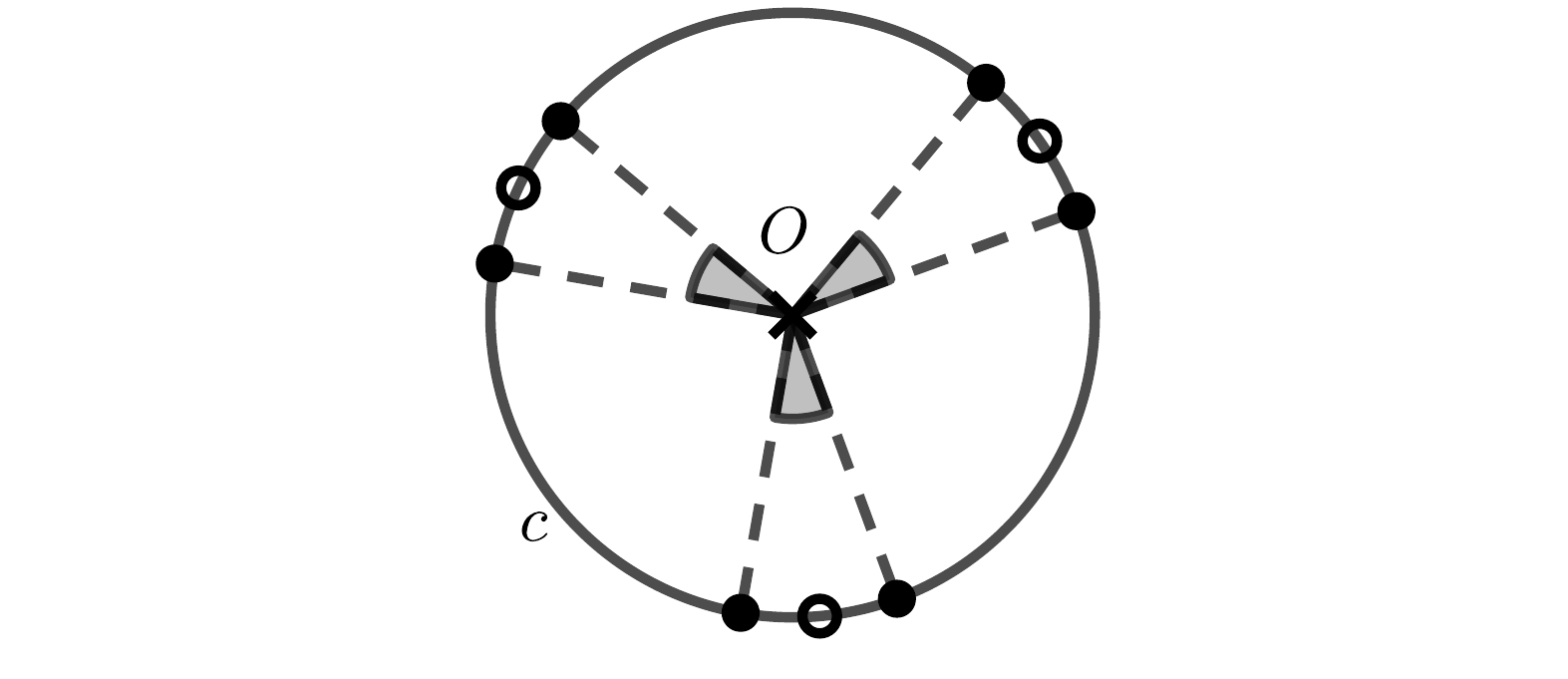}
\caption{Illustration of Algorithm 1. The set of black dots ($\bullet$) corresponds to the input and the set of white dots ($\circ$) to the output.}

\end{figure}

\noindent \begin{tabularx}{\textwidth}{|l X|} \hline  & \textbf{Algorithm 2 (definition of $\widetilde{\mathfrak A}_n$)} \\

\hline

INPUT:& a cyclic $n$-set   $\widetilde P\in\widetilde{\mathcal P}_n$ contained in a circle $c$ with circumcenter $O$, not containing any rotationally symmetric subset.\\

OUTPUT:& a point $\widetilde{\mathfrak A}_n(\widetilde P)$ in $c$.\\

\hline (1) & If $\widetilde P$ consists in a single point $V$, then $\widetilde{\mathfrak A}_n(\widetilde P)=V$.\\

(2) & If $\widetilde P$ consist in two points $V,W$, then $\widetilde{\mathfrak A}_n(\widetilde P)$ is the midpoint of the smallest circular arc among the two of them with endpoints $V,W$.\\

(3) & In other case, consider the pairs $\{V_1,W_1\},\ldots,\{V_r,W_r\}\subset \widetilde P$, consisting in two adjacent elements in $\widetilde P$ that maximize $\angle(V_i,O,W_i)$ (we take the angle corresponding to the circular arc not containing any other element in $\widetilde P$). \\

(4) & If $\widetilde P=\bigcup_{i=1}^r\{V_i,W_i\}$, then define a new $\widetilde P$ consisting in the midpoints of the circular arcs corresponding to the maximal angles described above and go back to Step (1).\\

(5) & In other case, define a new $\widetilde P$ obtained from the old one substracting  $V_1,\ldots,V_r,W_1,\ldots,W_r$ and go back to Step (1).\\

\hline

\end{tabularx}

\begin{figure}[h!]
\centering
\includegraphics[width=0.7\textwidth]{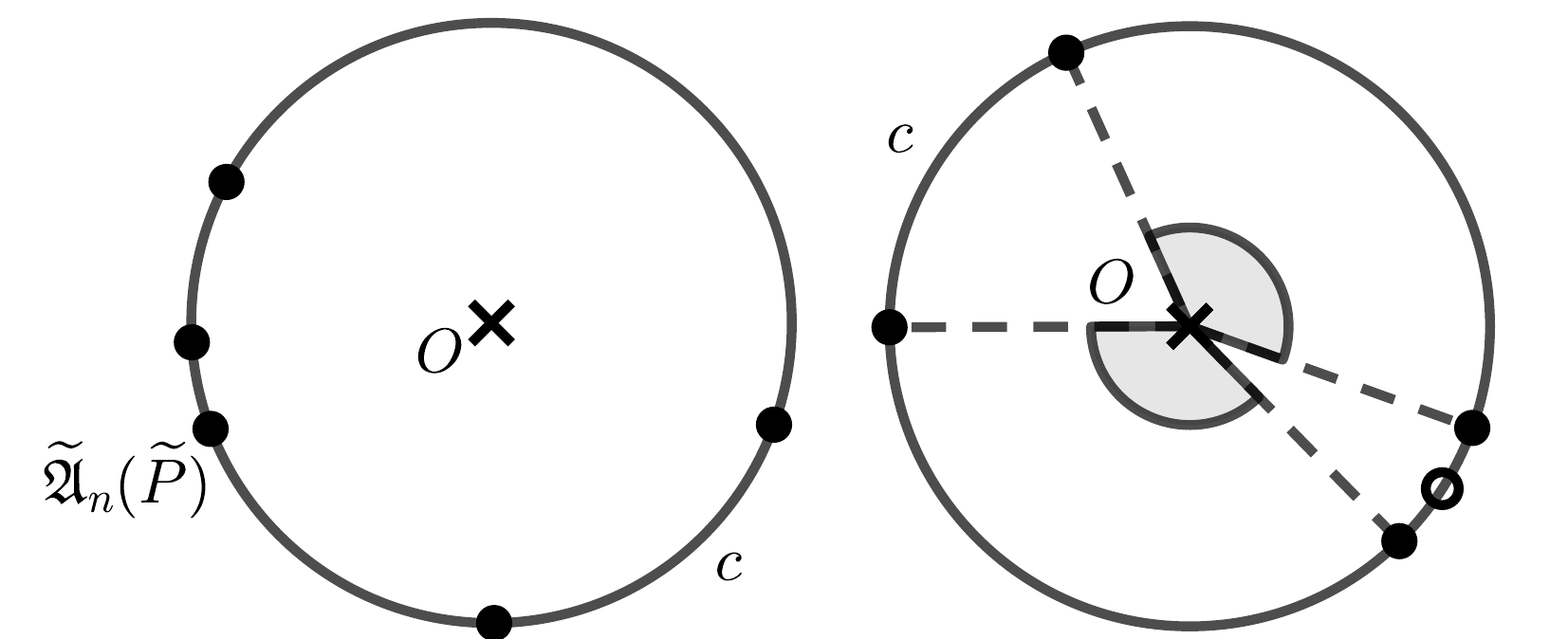}
\caption{Two examples corresponding to Algorithm 2, with the same convention as in the previous figure. In the example on the left, the point from the output coincides with one of the points in the set of the input.}

\end{figure}

\noindent \begin{tabularx}{\textwidth}{|l X|} \hline  & \textbf{Algorithm 3 (definition of $\widetilde{\mathfrak B}_n$)} \\ 
\hline INPUT:& a cyclic $n$-set $\widetilde P\in\widetilde{\mathcal P}_n$ with circumcenter $O$ and contained in the circle $c$.\\
OUTPUT:& a point $\widetilde{\mathfrak B}_n(P)$ in $c\cup \{O\}$ such that $\widetilde{\mathfrak B}_n(P)=O$ if and only if $\widetilde P$ is rotationally symmetric.\\

\hline (1)& If $\widetilde P$ is rotationally symmetric, then define $\widetilde{\mathfrak B}_n(\widetilde P)=O$.  \\
%(2)& In other case, let us  consider the sets $\widetilde Q_1,\ldots, \widetilde Q_n$ which are defined iteratively as follows. $\widetilde  Q_n$ contains the set of elements $V$ in $\widetilde P$ such that, for every $i\in\mathbb Z$,  the rotation of angle $\frac{i}{n}2\pi$ maps $V$ to another element in $ \widetilde P$. For $1\leq k\leq n-1$, $\widetilde Q_k$ contains the set of elements $V$ in $\widetilde P\setminus (\bigcup_{i=k+1^n}\widetilde Q_i)$ such that, for every $i\in\mathbb Z$,  the rotation of angle $\frac{i}{k}2\pi$ maps $V$ to another element in $\widetilde P\setminus (\bigcup_{i=k+1^n}\widetilde Q_i)$.\\

(2)& In other case, let us  consider the sets $\widetilde Q_2,\ldots, \widetilde Q_n$ such that  $2\leq i\leq n$, $\widetilde Q_i$ contains the set of elements $V\in\widetilde P$ such that, for every $i\in\mathbb Z$,  the rotation of angle $\frac{i}{k}2\pi$ maps $V$ to another element in $\widetilde P$.\\

(3)& Define $\widetilde Q_\infty=\widetilde P\setminus (\bigcup_{i=2}^n \widetilde Q_i)$.\\

(4)& If $\widetilde Q_{\infty}\neq \emptyset$, then $\widetilde{\mathfrak B}_n(\widetilde P)=\widetilde{\mathfrak A}_n(\widetilde P)$.\\

(5) & In other case consider, among all the pairs of integers $(i,j)$ such that $\widetilde Q_i,\widetilde Q_j\neq \emptyset$ and $gcd(i,j)=1$, the maximal one with respect to the lexicographic order.\\

(6) & Let $\widetilde Q=\Phi(\widetilde Q_i)\cup\Phi(\widetilde Q_j)$.\\

(7) & There are at most two pairs $\{V,W\}$ of adjacent points in $\widetilde Q$ minimizing the angle $\angle(V,O,W)$ (Statement (iii) in Lemma \ref{lemm.aux}). Suppose that this minimum is $\alpha$.\\

(8) & If there is only one such pair, then $\widetilde{\mathfrak B}_n(\widetilde P)$ is the midpoint of the circular arc with endpoints $V,W$ (the smallest one  among the two of them).\\

(9) & If there are two of them $\{V_1,W_1\}$ and $\{V_2,W_2\}$, then consider the two midpoints $M_1,M_2$ of the circular arcs with endpoints $V_1,W_1$ and $V_2,W_2$ (the one not containing any element in $\widetilde Q$).\\

(10)&  $M_1,M_2$ are not antipodal (Statement (iii) in Lemma \ref{lemm.aux}). So  $\widetilde{\mathfrak B}_n(\widetilde P)$ is the midpoint of the circular arc with endpoints $M_1,M_2$ (the smallest one among the two of them).\\

%(11) & If $M_1,M_2$ are antipodal, there are at most two pairs $\{V,W\}$ of adjacent points in $\widetilde Q$ minimizing the angle $\angle(V,O,W)$ among those greater than $\alpha$. Repeat Step (8), (9) and (10) (this time the midpoints cannot be antipodal, see Statement (ii) in Lemma \ref{lemm.aux}).\\

\hline

\end{tabularx}

\begin{figure}[h!]
\centering
\includegraphics[width=0.45\textwidth]{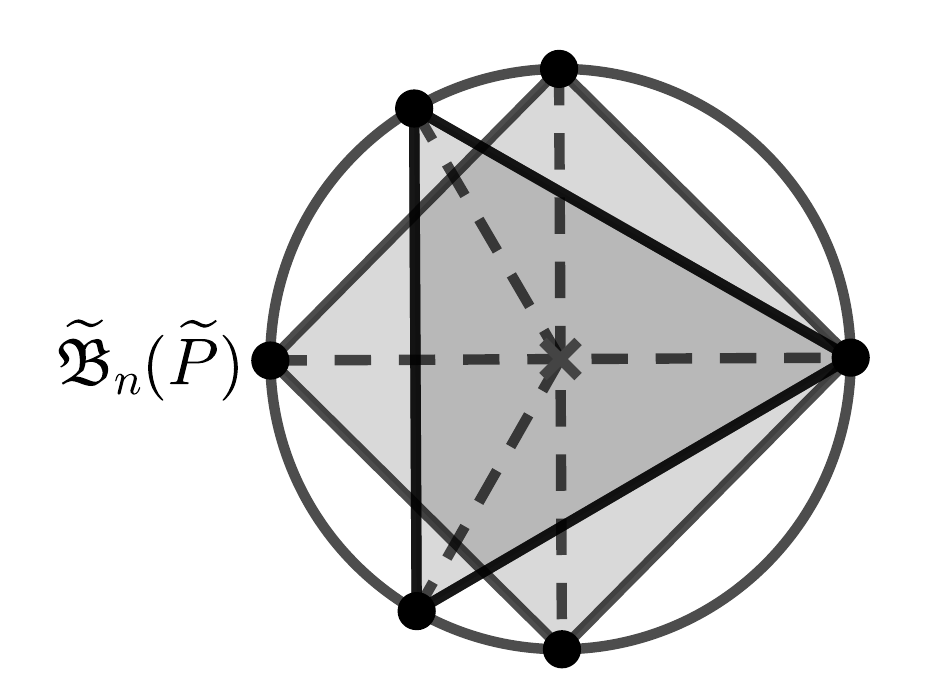}
\caption{Illustration of Algorithm 3. The output coincides with one of the points in the set in the input.}

\end{figure}

\noindent \begin{tabularx}{\textwidth}{|l X|} \hline  & \textbf{Algorithm 4 (definition of $\widetilde{\mathfrak B}_n^*$)} \\ 

\hline INPUT:& a ciclyc $n$-set with its vertices labelled $\widetilde P^*\in\widetilde{\mathcal P}_n^*$ with circumcenter $O$ and contained in the circle $c$. \\

OUTPUT:& a point $\widetilde{\mathfrak B}_n^*(\widetilde P^*)$ in $c\cup \{O\}$ such that $\widetilde{\mathfrak B}_n^*(\widetilde P^*)=O$ if and only if $\widetilde P^*$ is rotationally symmetric.\\

\hline  (1)& Suppose that all the elements in $\widetilde P^*$ have the same label. Then $\widetilde{\mathfrak B}_n^*(\widetilde P^*)=\widetilde{\mathfrak B}_n(\widetilde P)$, where $\widetilde P$ is the set obtained from $\widetilde P^*$ removing the labels.\\

(2)& Suppose that there are $k$ different labels, $k\geq 2$. Let us assume that they are $1,\ldots,k$ (in other case we do a relabelling preserving the order). Let $\widetilde Q_1^*,\ldots, \widetilde Q_k^*$ be the sets of points in $\widetilde P^*$ corresponding to the same label.\\

(3)& Suppose that at least one of the sets $\widetilde Q_1^*,\ldots,\widetilde Q_k^*$ is not rotationally symmetric, and suppose that $i$ is the maximal subindex with this property. Then $\widetilde{\mathfrak B}_n^*(\widetilde P^*)=\widetilde{\mathfrak B}_n(\widetilde Q_i)$, where $\widetilde Q_i$ is the set obtained from $\widetilde Q_i^*$ removing the labels.\\

(4) & In other case, for $i=1,\ldots,k$, define $d_i$ to be the order of the rotational symmetry.\\

(5)& Consider the pair of sets such that $gcd(d_i,d_j)=1$ and, corresponding to them, the maximal pair of subindices $(i,j)$ (with respect to the lexicographic order). \\

(6) & Let $\widetilde Q=\Phi(\widetilde Q_i)\cup\Phi(\widetilde Q_j)$, where $\widetilde Q_i,\widetilde Q_j$ are the sets obtained removing the labels. \\

& \\
 & Follow Steps (7), (8), (9), (10) and (11) in Algorithm 3 (replacing  $\widetilde{\mathfrak B}_n(\widetilde P)$ by $\widetilde{\mathfrak B}_n^*(\widetilde P^*)$).\\

\hline

\end{tabularx}

%%%%%%%%%%%%%%%%%%%%%%%%%%%%%%%%%%%%%%%%%%%%%%%%%%%%%%%%%%%%%%%
\section{Main result for $n$-multisets} \label{section.multiset}

In the following, let us denote by $\widetilde{\mathcal A}_n$, $\widetilde{\mathcal B}_n$ and $\widetilde{\mathcal C}_n$ to the elements in $\widetilde{\mathcal P}_n$ which group of symmetries has one fixed point, a line of fixed points and a plane of fixed points, respectively.

\begin{lemm}[definition of $\widetilde{\mathfrak X}_n$] \label{lemm.2} There exists a $n$-multiset center $\widetilde{\mathfrak X}_n:\widetilde{\mathcal P}_n\to\mathbb R^2$, that we will call the \textbf{center of rotational asymmetry}, such that $\widetilde{\mathfrak C}_n(\widetilde P)=\widetilde{\mathfrak X}_n(\widetilde P)$ if and only if $\widetilde P\in  \widetilde{\mathcal A}_n$.

\end{lemm}

\begin{proo} For every $\widetilde P\in\widetilde{\mathcal A}_n$, we define $\widetilde{\mathfrak X}_n(\widetilde P)=\widetilde{\mathfrak C}_n(\widetilde P)$. For each fixed $P=\{V_1,\ldots, V_n\}\in\widetilde{\mathcal B}_n\cup\widetilde{\mathcal C}_n$ we define $\widetilde{\mathfrak X}_n(\widetilde P)$ as follows.

 There is a finite set of  circles $c_1,\ldots,c_k$, labelled in decreasing order or radius (that can be 0) centered in $\widetilde{\mathfrak C}_n(\widetilde P)$ and intersecting $\widetilde P$.  There is a finite set of rays $r_1,\ldots, r_m$ with initial point $\widetilde{\mathfrak C}_n(\widetilde P)$ and intersecting $\widetilde P$.  

Define $X_1,\ldots, X_m$ where $X_i=c_1\cap r_i$. For each of these points $X_i$, consider the sequence $(a^i_1,\ldots, a^i_k)$ where $a^i_j$ is the number of points in $r_i\cap c_j$.

Define any set $\widetilde Q^*=\{(X_1,l_1),\ldots,(X_m,l_m)\}$ where $l_i$ is some natural number such that $l_i\leq l_j$ if and only if $(a^i_1,\ldots,a^i_k)\leq(a^j_1,\ldots,a^j_k)$ (in the lexicographic order). Then $\widetilde{\mathfrak X}_n(\widetilde P)=\widetilde{\mathfrak B}_n^*(\widetilde Q^*)$.

We can see that the function $\widetilde{\mathfrak X}_n:\widetilde{\mathcal P}_n\to\mathbb R^2$ defined according to the previous rules is a $n$-multiset center and $\widetilde{\mathfrak C}_n(\widetilde P)=\widetilde{\mathfrak X}_n(\widetilde P)$ if and only if $\widetilde P\in \widetilde{\mathcal A}_n$.

\end{proo}

\begin{lemm}[definition of $\widetilde{\mathfrak Y}_n$] \label{lemm.3}   There exists a $n$-multiset center $\widetilde{\mathfrak Y}_n:\widetilde{\mathcal P}_n\to\mathbb R^2$, that will be called \textbf{center of axial asymmetry}, such that $\widetilde{\mathfrak C}_n(\widetilde P)$, $\widetilde{\mathfrak X}_n(\widetilde P)$, $\widetilde{\mathfrak Y}_n(P)$ are pairwise distinct and not collinear if and only if $\widetilde P\in\widetilde{\mathcal C}_n$.

\end{lemm}

\begin{proo} For every $\widetilde P\notin \widetilde{\mathcal C}_n$, we define $\widetilde{\mathfrak Y}_n(\widetilde P)=\widetilde{\mathfrak C}_n(\widetilde P)$. For each fixed $\widetilde P=\{V_1,\ldots,V_n\}\in\widetilde{\mathcal P}_n$,  we define $\widetilde{\mathfrak Y}_n(\widetilde P)$ as follows.

Let  $r$ be the line passing through $\widetilde{\mathfrak C}_n(\widetilde P)$, $\widetilde{\mathfrak X}_n(\widetilde P)$.  The line $r$ divides the plane in two half-planes, that we will denote as $\mathcal H_1,\mathcal H_2$.

There is a finite set of concentric circles $c_1,\ldots, c_{k}$ labelled in decreasing order of radius (that can be $0$) centered in $\widetilde{\mathfrak C}_n(\widetilde P)$ and intersecting $P$.  There is a finite set of concentric circles $r_1,\ldots,r_{m}$ labelled in decreasing order of radius (that can be $0$) centered in $\widetilde{\mathfrak X}_n(\widetilde P)$ and intersecting $P$.

 Let us assign to each point $X\in \widetilde P$  the label $(i,j)$ if $X\in c_i\cap r_j$ and consider the two multisets of labels $L_1,L_2$, each of them corresponding to the labels of the points in $\mathcal H_1,\mathcal H_2$, respectively.

Let the total order relation $\leq_L$ defined as 
$$L_1\leq_L L_2\Longleftrightarrow min_{\text{lexic. order}}(L_1\Delta L_2)\in L_1. $$

\noindent  Define $v$ to be the free unit vector perpendicular to $r$ and pointing to the first half-plane if $L_1\leq_L L_2$ and the opposite one if $L_1\geq_LL_2$. Let ${\lambda=\sum_{V\in\widetilde P}\|V-\widetilde{\mathfrak C}_n(\widetilde P)\|}$. Then, we define $\widetilde{\mathfrak Y}_n(\widetilde P)=\widetilde{\mathfrak C}_n(\widetilde P)+\lambda v$.

The function $\widetilde{\mathfrak Y}_n:\widetilde{\mathcal P}_n\to\mathbb R^2$ defined according to the previous method is a $n$-multiset center and  $\widetilde{\mathfrak C}_n(\widetilde P)$, $\widetilde{\mathfrak X}_n(\widetilde P)$, $\widetilde{\mathfrak Y}_n(\widetilde P)$ are pairwise distinct and not collinear  if and only if $\widetilde P\in\widetilde{\mathcal C}_n$.

\end{proo}

%(1) & Let $\mathcal V$ be the multiset  consisting in the elements in $\mathcal U$  maximizing \\ & the distance from $\mathfrak C_n(P)$. \\ (2) & If $\mathcal V$ is symmetric with respect to $r$ define a new $\mathcal U$ obtained \\ & the old $\mathcal U$ substracting $\mathcal V$ and regurn to Step (1). \\ (3) & Let $\mathcal W$ be the multiset consisting in the elements in $\mathcal V$\\ & that maximize the distance from $r$. \\ (4) & If $\mathcal W$ is symmetric with respect to $r$, define a new $\mathcal V$ from the \\ & old one substracting the elements in $\mathcal V_n$ and return to Step (3) \\ 
%(5) & If there are more points of $\mathcal W$ in one of the two sides of $l$,\\ &  define the free unit vector $v$ as the one $v\perp l$ and pointing\\ & to this side. Jump to Step (9). \\ 
%(6) & In other case, the elements in $\mathcal W$ lie in the vertices of a rectangle.\\ & Suppose that the vertices are $X_1,X_2,X_3,X_4$ where \\& $X_1,X_3$ and $X_2,X_4$ lie in the same side of $r$\\ & and $X_1,X_2$ are further from $\mathfrak X_n(P)$ than $X_3,X_4$. \\ 
%(8) & Otherwise, there are a different number of points in $V_3$ and in $V_4$ \\ &  define the free unit vector $v$ as the one $v\perp l$ and pointing\\ & to this side.\\
%(9)& Define ${\lambda=\sum_{i=1}^n\|V_i-\mathfrak C_n(P)\|}$ and $\mathfrak Y_n(P)=\mathfrak C_n(P)+\lambda v$. \\\hline\end{tabular}

\begin{theo}[Main Problem for $n$-multisets] \label{theo.multiset} $ $

\begin{itemize}

\item[(i)] For every $n$-multiset center $\widetilde{\mathfrak Z}_n:\widetilde{\mathcal F}_n\to\mathbb R^2$ and every $\widetilde P\in\widetilde{\mathcal F}_n$, $\widetilde{\mathfrak Z}_n(\widetilde P)$ always belongs to the set of points fixed by the group of symmetries of $\widetilde P$. 

\item[(ii)] Conversely, for every point $X$ in the set of points fixed by the group of symmetries of a given $\widetilde P\in\widetilde{\mathcal P}_n$, there exists a $n$-multiset center $\widetilde{\mathfrak Z}_n:\widetilde{\mathcal P}_n\to\mathbb R^2$ such that $\widetilde{\mathfrak Z}_n(\widetilde P)=X$.
\end{itemize}

\end{theo}

\begin{proo} To prove the first statement, let $\widetilde{P}=\{V_1,\ldots, V_n\}\in\widetilde{\mathcal F}_n$, suppose that $T$ is in the symmetry group of $\widetilde P$. See that, by definition of center,
$$T(\widetilde{\mathfrak Z}_n(\widetilde P))=\widetilde{\mathfrak Z}_n(\{T(V_1),\ldots, T(V_n)\})=\widetilde{\mathfrak Z}_n(\{V_1,\ldots,V_n\})=\widetilde{\mathfrak Z}_n(\widetilde P).$$

For the second statement we need the fact that, for any set of centers $\widetilde{\mathfrak X}_n^1,\ldots,\widetilde{\mathfrak X}_n^k$ the affine combination
$$\lambda_1\widetilde{\mathfrak X}_n^1+\ldots+\lambda_k\widetilde{\mathfrak X}_n^k \qquad \text{ for }\lambda_1+\ldots+\lambda_n=1$$

\noindent is also a center. So, for each $\widetilde{P}\in\widetilde{\mathcal P}_n$, for every $X$ in the set of fixed points in its group of symmetries and for $\widetilde{\mathfrak C}_n$, $\widetilde{\mathfrak X}_n$, $\widetilde{\mathfrak Y}_n$, being the centroid, the center of rotational asymmetry and the center of axial asymmetry, there is at least one center
$$\widetilde{\mathfrak Z}_n=\lambda_1\widetilde{\mathfrak C}_n+\lambda_2\widetilde{\mathfrak X}_n+\lambda_3\widetilde{\mathfrak Y}_n \text{ for some }\lambda_1+\lambda_2+\lambda_3=1$$

\noindent such that $\widetilde{\mathfrak Z}_n(\widetilde{P})=X$.

\end{proo}
,

%%%%%%%%%%%%%%%%%%%%%%%%%%%%%%%%%%%%%%%%%%%%%%%%%%%%%%%%%%%%%%%%%%%%%%%%%
\section{Main result for $n$-gons} \label{section.polygon}

The main difficulty in this section is that we need to detect asymmetry in the adjacency relation, not in the vertices. In the following, let us denote by $\mathcal A_n$, $\mathcal B_n$ and $\mathcal C_n$ to the elements in $\mathcal P_n$ which group of symmetries has one fixed point, a line of fixed points and a plane of fixed points, respectively.

%For the two lemmas in this section, we need the following construction. A difficulty for us is that if a polygon $P$ is not simple, there is no natural notion of orientation for $P$.

%\begin{defi} Let $P\in\mathcal P_n$ such that the set of vertices $\widetilde P$ of $P$ is not collinear.  For any closed polygonal chain $(V_1,\ldots,V_{n+1})$, such that $V_1i\neq\mathfrak C_n(P)$ we say that it has a positive orientation if it satisfies the following:

%\begin{itemize}

%\item If there exists some $i$ such that $A_{n+1-i},A_1,A_{1+i}$ are not collinear, then for the  the first such $i$, the circular arc from $A_{n+1-i}$ to $A_{1+i}$ passing through $A_1$ is positively (resp. negatively) oriented.

%\item In other case, there exists some $i$ such that $A_{n+1-i},\ldots,A_{n+1},A_1,\ldots,A_{1+i}$ are not collinear and, for the first such $i$, the vector $\vec{V_2V_{n+1}}$

%circular arc with center $A_1$ from $A_2$ to $A_{1+i}$ is positively oriented.

%\end{itemize}

%\end{defi}

%Let $P\in\mathcal P_n$. For each vertex $W$ in $P$ there are two closed polygonal chains starting at $W$. If we denote the two different sequences as 
%$$({A_1},A_{2},\ldots,A_{n+1}),\qquad(B_1,B_{2},\ldots,B_{n+1}) $$

%\noindent then we have that, for $j=0,\ldots, n+1$, $A_{1+i}=B_{n+1-i}$. One of them is positively oriented and the other one negatively oriented. We will always denote by $({A_1},A_{2},\ldots,A_{n+1})$ to the positively oriented one and by  $(B_1,B_{2},\ldots,B_{n+1})$ to the negatively oriented one.

\begin{lemm}[definition of $\mathfrak X_n$] \label{lemm.2polygon} There exists a $n$-gon center ${\mathfrak X_n}:\mathcal P_n\to\mathbb R^2$, that we will call the \textbf{center of rotational asymmetry}, such that $\mathfrak C_n(P)=\mathfrak X_n(P)$ if and only if $P\in\mathcal A_n$.

\end{lemm}

\begin{proo} Let every $P\in\mathcal P_n$. Let us denote $\widetilde P$ to be the $n$-multiset of vertices of $P$. If $P\in\mathcal A_n$ we define $\mathfrak X_n(P)=\mathfrak C_n(P)$.  If $P\notin \mathcal A_n$ and  $\widetilde{P}\notin\widetilde{\mathcal A}_n$, then we define  $\mathfrak X_n(P)=\widetilde{\mathfrak X}_n(\widetilde{P})$. In other case, that is, $P\notin\mathcal A_n$ and $\widetilde P\in\widetilde{\mathcal A}_n$, we define $\mathfrak X_n(P)$ as follows.

The case in which $\widetilde P$ is collinear, requires a different approach. Let $\{W_1,\ldots,W_{m}\}$ be the multiset of elements at maximal distance from $\mathfrak C_n(P)$. Consider the set of of sequences of integers (2 for each element)
\begin{equation} \label{eq.seqcollinear}(a_{12},\ldots,a_{1,n+1}), (b_{12},\ldots,b_{1,n+1}),\qquad (a_{m2},\ldots,a_{m,n+1}), (b_{m2},\ldots,b_{m,n+1}),\end{equation}

\noindent each pair corresponding to the two polygonal sequences $(A_{i1},\ldots,A_{i,n+1})$, $(B_{i1},\ldots,B_{i,n+1})$ starting at each $W_i$ and such that $a_{ij}$ (resp. $,b_{ij}$) denote the distance from $A_{i,j-1}$ to $A_{ij}$ (resp. from $B_{i,j-1}$ to $B_{ij}$) with positive sign if $A_{ij}$ is closer to $\mathfrak C_n(P)$ than $A_{i,j-1}$. These sequences in Equation \eqref{eq.seqcollinear} are all of them different. So  only one of these sequences is minimal with respect to the lexicographic order. The point $W_i$ corresponding to these sequence will be $\mathfrak X_n(P)$.

Suppose that not all the points in $\widetilde P$ are collinear. Define $c_1,\ldots,c_k$ and $r_1,\ldots,r_m$ as in the proof of Lemma \ref{lemm.2}.

Let $\{W_1,\ldots,W_{m_1}\}$ be the multiset which elements are in $P\cap c_1$.  For each $W_i$, consider the two polygonal sequences starting at $W_i$ and denote them by
$$(A_{i1},\ldots,A_{i,n+1}),\qquad(B_{i1},\ldots,B_{i,n+1}). $$

If the polygon $P$ is not simple, we do not have a natural notion of orientation. But in this setting, we will be able to  stablish a criteria to say that one of these two sequences is positively oriented and the other one negatively oriented.  By hypothesis $\widetilde P$ is not collinear. Consider $r$ to be the ray from $\mathfrak C_n(P)$ to $W_i$. Let $r_+$ (resp. $r_-$) to be the ray from $\mathfrak C_n(P)$ to some of the points in $P\setminus r$ such that the angle from the ray $r$ to $r_+$ goes in the positive (resp. negative) sense and is minimal. From all the points of $\widetilde P$ in $r_+$ (resp. $r_-$), consider the furthest one from $\mathfrak C_n(P)$ and denote it by $V_+$ (resp. $V_-$). In the following, we suppose that $(A_{i1},\ldots,A_{i,n+1})$ reaches $V_+$ before $V_-$ and we will call it the positively oriented sequence. $(B_{i1},\ldots,B_{i,n+1})$ will reach $V_-$ before $V_+$ and will be called the negative oriented sequence.

\noindent Associated to this positively oriented sequence, we are going to define another sequence $(a_{i2},\ldots,a_{i,n+1})$.  For $j=2,\ldots,n$, the element $a_{ij}$    is a pair  $(\rho,\alpha)$ where $A_{ij}\in c_\rho$ and $\alpha$ is the angle $\angle(A_{ih},\mathfrak C_n(P),A_{ij}) $ where $h<j$ is the last index such that $A_{ih}\neq\mathfrak C_n(P)$ and we consider that $\angle(A_{ih},\mathfrak C_n(P),A_{ij})=0$ if $A_{ij}=\mathfrak C_n(P)$. We define a sequence. $(b_{i2},\ldots,b_{i,n+1})$ for the negatively oriented sequence in a similiar way with the corresponding modifications. See Figure \ref{P1}.

\begin{figure}[h!]
\centering

\includegraphics[width=0.5\textwidth]{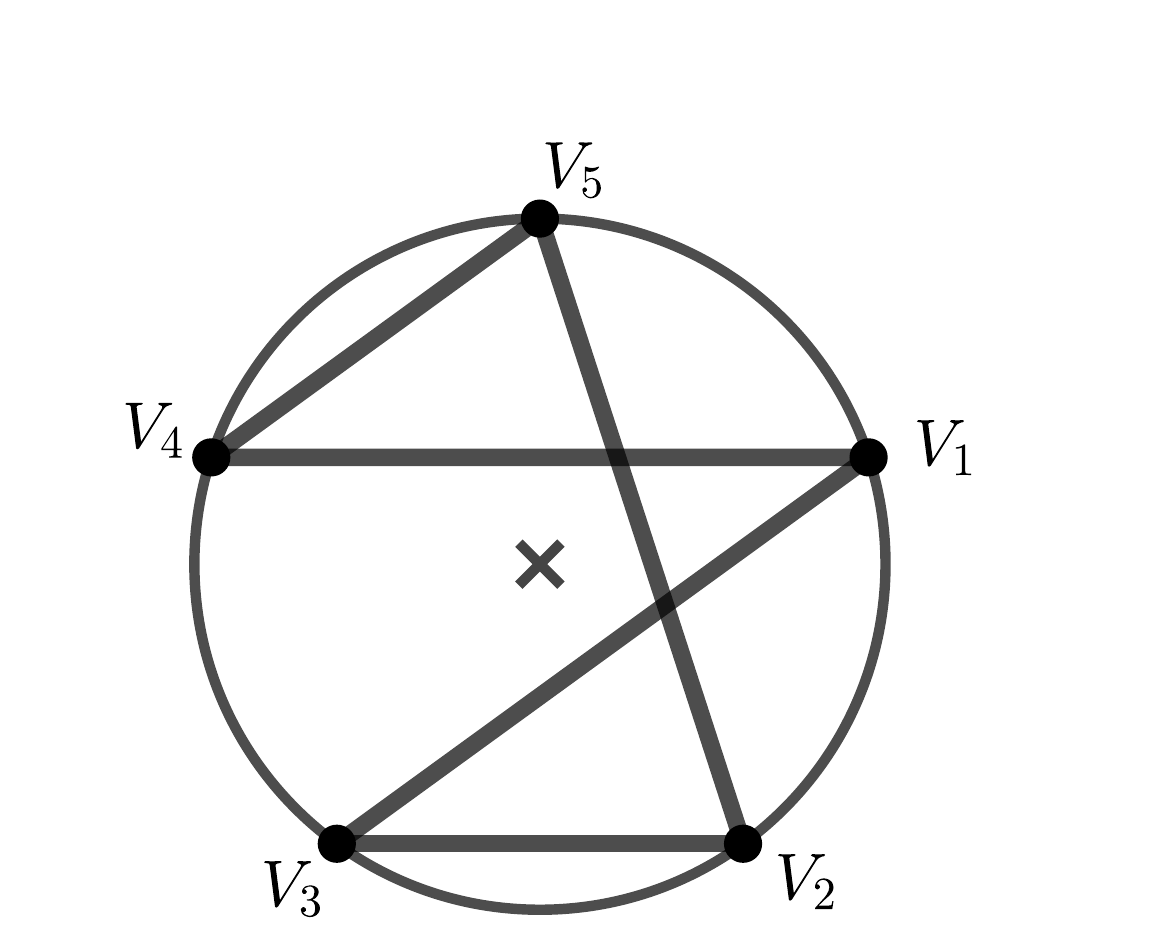}

\caption{In this example, there is a unique circle $c_1$. The positively oriented polygonal chain starting at $V_1$ is $(V_1,V_3,V_2,V_5,V_4,V_1)$ (in this case, $V_+=V_2$) and the corresponding sequence of pairs is $((1,2\frac{2\pi}{5}),(1,4\frac{2\pi}{5}),(1,3\frac{2\pi}{5}),(1,4\frac{2\pi}{5}),(1,2\frac{2\pi}{5}))$.}
\label{P1}

\end{figure}

Note that a $n$-gon is rotationally symmetric if and only if, for any $i$, the corresponding $(a_{i2},\ldots,a_{i,n+1})$ is periodic. As a consequence and since $P$ is not rotationally symmetric, then all the sequences must be different. Pick the point $Z_+\in\{W_1,\ldots,W_{m_1}\}$ which sequence is minimal with respect to the lexicographic order (using, in turn, the lexicographic order to compare each position). Define $Z_-$ in a similar way.

Define $\mathfrak X_n(P)$ as the midpoint of the circular arc with endpoints $Z_+,Z_-$. We choose the one of the two possible such that from $Z_+$ to $Z_-$ is positively oriented (and from $Z_-$ to $Z_+$ is negatively oriented).

\end{proo}

\begin{lemm}[definition of $\mathfrak Y_n$] \label{lemm.3polygon}  There exists a $n$-gon center ${\mathfrak Y_n}:\mathcal P_n\to\mathbb R^2$, that we will call the \textbf{center of rotational assymetry}, such that $\mathfrak C_n(P)$, ${\mathfrak X_n}(P)$, ${\mathfrak Y_n}(P)$ are different and not collinear if and only if $P\in\mathcal C_n$.

%Consider the set $\widetilde{\mathcal C_n}\cup\mathcal B_n$ and $\mathfrak X_n$ be the $n$-multiset center (and so a $n$-gon center) described in the corresponding lemma. Then there exists a $n$-gon center $\widetilde{\mathfrak Y_n}:\widetilde{\mathcal C_n}\cup\mathcal B_n\to\mathbb R^2$ such that, for every $P\in\widetilde{\mathcal C_n}\cup\mathcal B_n$, $\mathcal C_n(P),\mathfrak X_n(P),\widetilde{\mathfrak Y_n}(P)$ are not collinear.

\end{lemm}

\begin{proo} Let $P\in\mathcal P_n$ and let $\widetilde P$ denote the $n$-multiset of its vertices.

If $P\notin\mathcal C_n$, we define $\mathfrak Y_n(P)=\mathfrak C_n(P)$. If $P\in\mathcal C_n$ and $\widetilde P\in\widetilde{\mathcal C}_n$, then $\mathfrak Y_n(P)=\widetilde{\mathfrak Y}_n(\widetilde{P})$. In other case, that is, $P\in\mathcal C_n$ and $\widetilde P\notin\widetilde{\mathcal C}_n$, then we define $\mathfrak Y_n(P)$ as follows.

Let us define $r$, $\mathcal H_1,\mathcal H_2$, $c_1,\ldots,c_k$, $r_1,\ldots,r_m$ as in the proof of Lemma \ref{lemm.3}, using $\mathfrak C_n(P)$, $\mathfrak X_n(P)$ (instead of $\widetilde{\mathfrak C}_n(\widetilde P)$, $\widetilde{\mathfrak X}_n(\widetilde P))$. %Let us define, additionally, the points $W_{ij}=c_{i}\cap r_j$.

Let $Q=\{W_1,\ldots,W_l\}$ to be the elements in $c_1\cap r_1$ ($Q$ is a multiset with at most two different elements, both of them out of $r$).  For each $W_i\in Q$, consider the polygonal chains with root at $W_i$
$$(A_{i1},\ldots,A_{i,n+1}),\qquad(B_{i1},\ldots,B_{i,n+1}) $$

%$(2\leq j\leq\lfloor\frac{n-1}{2}\rfloor$ consider $A_{i,j+1},B_{i,j+1}$ to be the points adjacent to $A_{ij},B_{ij}$ respectively, different from $A_{i,j-1},B_{i,j-1}$. The sequences are polygonal chains. We say that $A_{ij_1}$ is a predecessor of $A_{ij_2}$ and that $A_{ij_2}$ is a sucessor of $A_{i,j_1}$ if $i_1<i_2$ (similar for $B_{ij_1},B_{ij_2}$).

Now define, associated to  each of these elements $W_i$ a couple of sequences $(a_{i2},\ldots,a_{i,n+1})$, $(b_{i2},\ldots,b_{i,n+1})$.  For $j=1,\ldots,n$, the element $a_{ij}$  (similar for $b_{ij}$) is a triple  $(x,y,z)$ such that $A_{ij}\in c_x\cap r_y$ and $z$ is defined as follows.
$$z=\begin{cases} -1& \text{if the last predecessor of $A_{ij}$ out of $r$ is not in the same halfplane}\\ & \text{as $A_{ij}$ or the earliest sucessor of $A_{ij}$ out of $r$}\\
 +1& \text{if the last predecessor of $A_{ij}$ out of $r$ is  in the same halfplane}\\ & \text{as $A_{ij}$ or the earliest sucessor of $A_{ij}$ out of $r$}\end{cases} $$

%There is a total order $\leq S$ to compare such sequences: the lexicographic order where each position $a_{ij}$ or $b_{ij}$ is compared using, in turn, the lexicographic order for the triples $(x,y,z)$. 

Note that $P$ is axially symmetric if for any (for all) $i=1,\ldots, k$ we have that $(a_{i2},\ldots,a_{i,n+1})=(b_{i2},\ldots,b_{i,n+1})$. Let us assume that, for every $i=1,\ldots,k$, $(a_{i2},\ldots,a_{i,n+1})$ is smaller than $(b_{i2},\ldots,b_{i,n+1})$ in the lexicographic order where each position $a_{ij}$ or $b_{ij}$ is compared using, in turn, the lexicographic order for the triples $(x,y,z)$.

Moreover, for al $i=1,\ldots,k$ the sequences $(a_{i2},\ldots,a_{i,n+1})$ are all different. Let us choose $W_i$ to be the point in $Q$ corresponding to the smallest sequence with respect to the order described above. Then, we define $\widetilde{\mathfrak Y}_n(\widetilde P)=W_i$.

The function ${\mathfrak Y}_n:{\mathcal P}_n\to\mathbb R^2$ defined according to the previous method is a $n$-gon center and  ${\mathfrak C}_n( P)$, ${\mathfrak X}_n(P)$, ${\mathfrak Y}_n(P)$ are pairwise distinct and not collinear  if and only if $\widetilde P\in{\mathcal C}_n$.

\end{proo}

%(1) & Let $\mathcal V$ be the multiset  consisting in the elements in $\mathcal U$  maximizing \\ & the distance from $\mathfrak C_n(P)$. \\ (2) & If $\mathcal V$ is symmetric with respect to $r$ define a new $\mathcal U$ obtained \\ & the old $\mathcal U$ substracting $\mathcal V$ and regurn to Step (1). \\ (3) & Let $\mathcal W$ be the multiset consisting in the elements in $\mathcal V$\\ & that maximize the distance from $r$. \\ (4) & If $\mathcal W$ is symmetric with respect to $r$, define a new $\mathcal V$ from the \\ & old one substracting the elements in $\mathcal V_n$ and return to Step (3) \\ 
%(5) & If there are more points of $\mathcal W$ in one of the two sides of $l$,\\ &  define the free unit vector $v$ as the one $v\perp l$ and pointing\\ & to this side. Jump to Step (9). \\ 
%(6) & In other case, the elements in $\mathcal W$ lie in the vertices of a rectangle.\\ & Suppose that the vertices are $X_1,X_2,X_3,X_4$ where \\& $X_1,X_3$ and $X_2,X_4$ lie in the same side of $r$\\ & and $X_1,X_2$ are further from $\mathfrak X_n(P)$ than $X_3,X_4$. \\ 
%(8) & Otherwise, there are a different number of points in $V_3$ and in $V_4$ \\ &  define the free unit vector $v$ as the one $v\perp l$ and pointing\\ & to this side.\\
%(9)& Define ${\lambda=\sum_{i=1}^n\|V_i-\mathfrak C_n(P)\|}$ and $\mathfrak Y_n(P)=\mathfrak C_n(P)+\lambda v$. \\\hline\end{tabular}

\begin{theo}[Main Problem for $n$-gons] \label{theo.polygon}$ $

\begin{itemize}

\item[(i)] For every $n$-gon center ${\mathfrak Z}_n:{\mathcal F}_n\to\mathbb R^2$ and every $P\in{\mathcal F}_n$, ${\mathfrak Z}_n(P)$ always belongs to the set of points fixed by the group of symmetries of $P$. 

\item[(ii)] Conversely, for every point $X$ in the set of points fixed by the group of symmetries of a given $P\in{\mathcal P}_n$, there exists a $n$-gon center ${\mathfrak Z}_n:{\mathcal P}_n\to\mathbb R^2$ such that ${\mathfrak Z}_n(P)=X$.
\end{itemize}

\end{theo}

\begin{proo} The proof is exactly the same as in Theorem \ref{theo.multiset}, but removing the tildes.
\end{proo}
,

\section{Final comments}

The intention of this notes is to provide a theoretical result relating the group of symmetries of a given multiset or polygon and the location of the centers of this object.

The centers $\widetilde{\mathfrak X}_n$, $\widetilde{\mathfrak Y}_n$, $\mathfrak X_n$, $\mathfrak Y_n$ are defined in terms of the center $\widetilde{\mathfrak B}_n$. This last multiset center is not continuous, with respect to the most natural candidates for a topology for $\widetilde{\mathcal P}_n$, $\mathcal P_n$, (such as the one induced by the Hausdorff distance for $\widetilde{\mathcal P}_n$ and the corresponding quotient topology for $\mathcal P_n$). This make that the centers $\widetilde{\mathfrak X}_n$, $\widetilde{\mathfrak Y}_n$, $\mathfrak X_n$, $\mathfrak Y_n$ look rather artificial.

Now that the main results in this paper are proved and we know the answer for our Main Problem,  we could look for sets of centers $\widetilde{\mathfrak X}^1,\ldots,\widetilde{\mathfrak X}^n$ and $\mathfrak X^1,\ldots,\mathfrak X^n$, among the ones with geometric meaning, such that if all of them coincide (resp. are collinear) the corresponding object has rotational symmetry (resp. axial symmetry), in a more similar fashion to the  ``inspiring results'' listed in the introduction. To do so, the difficulty lies in solving equations of the type
$$\widetilde{\mathfrak Z}^1(\widetilde P)=\widetilde{\mathfrak Z}^2(\widetilde P) \qquad \text{or}\qquad \mathfrak Z^1(P)=\mathfrak Z^2(P)  $$

\noindent in the indeterminate $\widetilde P$, $P$, respectively. This equation require some algebraic study and will be a matter for future work.

%\begin{acknowledgements}
%\end{acknowledgements}




\end{document}